# Linear programming based evacuation models for a controlled freeway

Sheng-Xue He[*]

*Business School, University of Shanghai for Science and Technology, Shanghai 200093, China*

**Abstract**: A linear programming (LP) model is proposed to improve the performance of a controlled freeway during an emergency evacuation. Based on reasonable assumptions, the main relationships among key factors are kept without the uncertain impact of subordinate factors in this model. Three vital issues related to optimal freeway control during an emergency evacuation are addressed. How to properly realize some predetermined priorities assigned to some on-ramps during an emergency evacuation is analyzed through sensitivity analysis of objective function coefficients. The time-varying throughput capacity of freeway mainline and the storage capacities of on-ramps are also addressed by using sensitivity analysis of the right-hand side (RHS) of corresponding constraints. The last key issue we focused on is how to deal with the uncertain demand in a crude form by applying robust optimization. The methodology presented in this paper can help evacuation managers realize a coordinated ramp metering with predetermined priorities put on some on-ramps during an emergency evacuation. The result obtained by analyzing the impact of time-varying mainline capacity provides the evacuation manager a chance of using the extra capacity in some segments of mainline to buffer the traffic so as to reduce the congestion related to upstream on-ramps. To reduce the risk due to the possible infeasibility of deterministic control approaches in an uncertain situation, we apply an affinely adjustable robust counterpart (AARC) approach to deal with uncertain dynamic demands restricted to polyhedral sets. The preliminary numerical experiments show that the AARC approach provides us with a promising solution comparing to the sampling based stochastic approach.

## 1. Introduction

Increasing disasters experienced all around the world during the past several decades have made the emergency evacuation an important issue for both practitioners and theorists. If an evacuation needs to be deployed in a region within a pressing time period, the traffic demand generated will generally exceed the available network capacity. So congestion, even serious incidents, may arise at this time due to substantial traffic demand. Freeway, used as the main infrastructure during regional emergency evacuation, supplies the main channel or backbone in the whole evacuation network. Though there are many researches about freeway control and management in the existing literature, the number of those specially related to a controlled freeway under an emergency evacuation situation is relative small (So and Daganzo, 2010). Various control methods have been studied under the ordinary traffic situation. As is well known, ramp metering is one of the most important methods. The others include variable speed limits, variable message signs (route guidance), contra-flow lane, et al. A coordinated and integrated control strategy is more difficult to be designed and implemented than a strategy using only one individual control method. Though the knowledge about freeway control known to us is substantial nowadays, a gap still exists between researches focused on an ordinary circumstance and those under an emergency evacuation situation.

We notice that the generic ramp metering or control strategies are often unfitted in an emergency evacuation. When considering a regional evacuation using a freeway or freeways' network, we assume that a completely controlled ramp system aiming at system optimization is generally required. To coordinate and integrate the available various control methods is desired to cope with the generally hard predictable and large demand. Proper evacuation preparedness should be considered to account for an emergency evacuation. For example, proper predetermined priorities of releasing flow should be assigned to some ramps with HOV or exclusive bus lanes. In view of the complexity and urgency of regional evacuation by a freeway, a simplified and efficient formulation of our problem using linear programming would be promising and desired.

---
[*] Tel.: +086-21-18019495040.
  E-mail address: lovellhe@126.com

Two observations confirm us the possibility of formulating our problem in a simple LP model. First, most of existing researches have shown that to prevent the internal queue in a freeway and to keep the mainline flow under the critical capacity will bring us the highest performance of a freeway (Wattleworth, 1967; Papageorgiou, 1998; Kotsialos, et al. 2002; Kotsialos and Papageorgiou, 2004; So and Daganzo, 2010). Based on this observation, many effective ramp metering strategies have been proposed and tested. Second, when considering the real life situation in an emergency evacuation, we find out it is readily satisfied that a freeway can be nearly completely controlled by ramp-metering though the uncertain large demand and the predetermined priorities to some specified ramps exist. If the second observation is true, then the first could be realized. With both, a linear programming mathematical model can be formulated assuming that a free flow status keeps on in the freeway mainline during an emergency evacuation. As is well known, the advantage of linear programming relies on its substantial existing research results and many powerful algorithms available to solve an LP with a large scale.

In our study, we will focus on three very important questions about freeway evacuation. They can be investigated taking advantage of the existing substantial knowledge about LP. The first one is how to efficiently use the available on-ramps' storage capacities and mainline reserved capacities. Intuitively we could expect to solve this problem in a completely controlled freeway by coordinating ramp metering and making use of the residual capacity of mainline to buffer more vehicles so as to mitigate the burden of queuing in on-ramps. The second one is about realizing the predetermined on-ramp discharging priority. If there are some special requirements that give the priorities of discharging flow to some on-ramps, e.g. bus only or HOV on-ramps, how should we deal with this situation to attain a system optimal performance? This problem could be addressed elegantly using sensitivity analysis of LP. The last important issue is to decide the optimal ramp metering rates under a situation with only uncertain demand information in a crude form. Generally speaking, the uncertain traffic demand in an emergency evacuation is time-dependent and could be restricted to some box or polyhedral constraints. How to design a ramp metering strategy to cope with this uncertainty is a challenge not only to practitioners but also to theorists. Recently, robust optimization of LP has witnessed a significant growth. The related studies provide a strong basis for our investigation in the last problem with which we are seriously concerned.

The history of using LP in freeway control field can be traced back to Wattleworth (1967). Though nowadays some more complex high order macroscopic traffic flow models are generally adopted by researchers in this field, one famous linear difference equation model, i.e. Cell Transmission Model (CTM), is also broadly used (Daganzo, 1994 and 1995). Based on CTM, a series of LP models related to freeway and urban networks have been proposed. To construct an LP model of freeway control, some special assumptions will needed. Either fixed travel time in the freeway mainline or slowly changing average speed is always assumed. Sometimes, stable and smooth dispersing rates related to off-ramps are also supposed (Zhang and Levinson, 2004). Temporal-spatial discrete difference equations instead of differential equations are used in some researches resulting from the adoption of Piecewise linear approximation of high order fundamental relationships between traffic flow and density. Compared to the existing studies, the models considered in this paper are closely related to the emergency evacuation situation and more simplicity and efficiency with respect to our main concerns. A brief clarification about the difference will be given in Section 2.

The main contributions of this paper are as follows. (a) LP optimization models formulated emphasize the main concerns and relationships under an emergency evacuation situation. With several reasonable assumptions, the uncertain influence of the other subordinate factors could be removed from the final formulation. Using the range analysis of the LP objective function coefficients to determine the proper discharging order and metering rates of on-ramps can realize the predetermined priorities given to some specified on-ramps. (b) Another benefit resulted from the concise LP model is the possibility of making full use of the available capacity of freeway. In the planning or preparedness stage, through solving the LP model, we can identify the potential bottlenecks in the freeway, and then take proper measures to effectively improve the related capacities. During the evacuation, the capacities of some mainline stretches could be made full use of to buffer more vehicles. To achieve this aim, we need first to use our model to identify the mainline segments with residual capacity, and then compute the possible buffer quantity. Through adjusting the related weights on waiting time in on-ramps, we expect to hold more vehicles in some stretches of freeway mainline where evacuees will face lower risk comparing to waiting at some upstream on-ramps. The time-variant capacity due to some control methods, e.g. contra-flow lane, can be treated easily using our LP model. Coordinated ramp metering strategy will be used to improve the performance of whole system with regard to the storage capacities of on-ramps. (c) Uncertainty of traffic demand can be coped with by robust optimization of LP elegantly. An optimal solution to the uncertain situation will be presented by our analysis to avoid the possibility of non feasible solution resulted from conventional methods. All above mentioned can help evacuation managers effectively and efficiently to manage a freeway in an emergency regional evacuation. Note that it is possible to combine with other existing ramp-metering methods to implement the results obtained in our study. For example, Advanced Motorway Optimal Control (AMOC) or CTM based freeway control system may be used with input from our model.

The paper is organized as follows. In the next section, a brief literature review will be presented. Section 3 will be devoted to formulating the freeway evacuation system in LP mathematical models. Section 4 will use the sensitivity analysis of LP to address issues of the discharging priority and the full use of throughput capacity. In Section 5, the uncertainty of demand in a hard form will be treated by using robust optimization. The numerical analysis will be carried out in Section 6. In Section 7, we will summarize the paper and give some recommendations for future research.

## 2. Literature review

In view of the topic of completely controlled freeway during an emergency evacuation, the related literature is limited. But about the optimal control of freeway, there are so many existing. Su and Daganzo (2010) propose an adaptive strategy for managing a freeway in an emergency evacuation. Their strategy gives priorities of releasing vehicles to upstream on-ramps. Through assuming a specified relation between the capacity of an on-ramp and the one of its immediately downstream mainline segment, they prove that the shortest clearance time can be obtained using their strategy. To implement their strategy, the quantity of dynamic demands along the freeway is unnecessary. Several assumptions used in this literature will be adopted in our study to make the subsequent analysis easy and concise, e.g. the infinite free flow speed of the mainstream. Zhang and Levinson (2004) have also attempted to optimally control a freeway without the origin-destination demand information. In their study, the stable dispersing rates of off-ramps are assumed. Though at first glance the method of Zhang and Levinson (2004) and the approach of So and Daganzo (2010) seem total different, a further careful observation will show there are many similarities between them. For example, two stages of implementation of their strategies are similar. From the latter, we take as truth their announcement that no effective ways existing could provide us with the accurate dynamic traffic demands. So a polyhedral type of demand constraints in a crude form will naturally be considered in our study.

To consider the macroscopic traffic flow models, an effective formulation to control freeway traffic is widely accepted (Kotsialos, et al. 2002; Kotsialos and Papageorgiou, 2004). This method evolves over time, from simple local control of a single ramp to complex coordinated and integrated control of a whole freeway network. Though we will not use this specific traffic flow model, several very important observations from their researches will be accepted and adopted by authors of this paper to form our basic assumptions. For example, the best performance could always be expected to happen under the condition that the mainline capacity does not be exceeded. The limitation of the storage capacity of on-ramp is a key factor in these models, which we will adopt in our formulation.

Linear programming technique was used in freeway control field by Wattleworth (1967). In this early study, an elegant model of peak-period freeway control is constructed under the assumption of static traffic input. Nowadays, a popular LP formulation, i.e. CTM, is widely used in traffic flow studies. Ziliaskopoulos (2000) proposed an LP model to deal with dynamic traffic assignment with system optimization as its objective. In his study, only the problem with a single destination was treated. Later many studies along this direction were conducted (e.g. Lo, 1999; Szeto and Lo, 2004; Gomes and Horowitz, 2006). Using CTM, some researchers have addressed issues related to emergency network evacuation (Chiu et al., 2007). Though there are many benefits to adopt CTM, we still feel it is too complicated to be adopted to deal with our concerns at the first stage of study in view of the circumstance of emergency evacuation. But we can see it is possible to take CTM as the second step to further our study in the future.

During an emergency network evacuation, some roads may have extra capacities though some may have used up their capacities and become bottlenecks of the evacuation network. These residual capacities in terms of extra available space in some time periods make it possible to lighten the upstream congestion or waiting time by holding more vehicles in these roads with extra capacities left unused. Generally speaking, the place nearer to the risk origin has higher potential risk level. So if we could move more evacuees to downstream roads with lower risk levels without significant impact on the performance of the whole system, the potential damage and loss due to the relative high risk level near the risk origin would be decreased to a relative low level. Duanmu et al. (2012) has addressed this issue. In their study, Duanmu et al. (2012) present an algorithm that allows a long link to be used as a buffer to keep the traffic flow moving in. Their aim is to keep the traffic under saturated in the buffer zone and minimize the total travel time simultaneously. In our study, we focus on how to identify the location where we may hold more vehicles. A primitive approach is also proposed to estimate the quantity of vehicles which may be buffered using the fundamental relationship between traffic flow and density.

To deal with the uncertain dynamic traffic demand during an emergency evacuation, we introduce the affinely adjustable robust counterpart (AARC) approach. This approach has been proposed by Ben-Tal et al. (2004) to address the "wait and see" type of decisions. For a summary of the state of art of robust optimization, readers can refer to Ben-Tal et al. (2009), Bertsimas et al. (2007) and references therein. The related techniques of robust optimization have been proposed to solve the network and transportation problems (e.g. Bertsimas and Perakis, 2005; Ordonez and Zhao, 2007; Atamturk and Zhang, 2007; Mudchanatongsuk et al., 2008; Erera et al., 2009; Yin et al., 2008, 2009). Recently, a series of researches using robust optimization with CTM as their modeling basis have been conducted by Yao et al. (2009), Ben-Tal et al. (2011) and Chung et al. (2012). Regarding to the difficulty of predicting and estimating the accurate dynamic traffic demand for an emergency evacuation, we think robust optimization, especially AARC, will become an effective and efficient approach to be used more widely in the future.

## 3. Model formulation

In this section, we will first introduce the basic assumptions adopted in our study and clarify the necessity of using them under an emergency evacuation situation. Then after the main notations are explained in subsection 3.2, we propose two basic models with different objectives in subsection 3.3. In subsection 3.4, we extend the basic models to cope with more realistic situations with new constraints added.

*3.1. The basic assumptions*

There are several assumptions that are special for freeway emergency evacuation. After the introduction of these assumptions, the difference between our control strategy and the ordinary ramp metering strategies will be clear.

First, we assume that the freeway to be used is empty when the evacuation begins. This assumption can be satisfied in real life when the evacuation begins in the early morning or at some time when the traffic is light in the freeway. This assumption will make the modeling easier and result in insignificant impact on the following analysis, simultaneously.

Second, we assume that the free flow speed is infinite in the freeway. The influence of this simplification is trivial on the final results according to So and Daganzo (2010). The benefit from this assumption is to make subsequent modeling and analyses concise and easy and to obtain some clear and convincible insights into freeway evacuation.

In real life, the changing vehicle speed is hard to be described precisely even with some microscopic models which are generally time-consumed and application-limited. The varying speed will restrict the accuracy of a model with fixed travel speed as one of its assumptions. Note that it is possible to extend the model assuming infinite free flow speed to one dealing with fixed speed. The extension in this direction is straightforward according to the current state of art. One way to realizing this extension is to synchronize the starting time along the freeway (Zhang and Levinson, 2004). This method is to make the time counts of different on-ramps start from different time instants, so as to realize that the released flow with same time label according to their starting time point will meet and merge on the mainline. After such synchronization, our model can be used to describe the freeway system with fixed travel speed as well.

Considering the assumption from a different angle, we will find out that the fixed finite speed and the infinite free flow speed are in fact equivalent in terms of impacts on the objective value. Supposing that the travel time on any link is fixed, we will find out that in the objective function the Total Evacuation Travel Time (TETT) is meaningfully determined by the Total Waiting Time (TWT) at on-ramps. The reason lies in that the total travel time after the vehicles enter the controlled freeway mainline is a constant for the clearance of network. This fixed travel time in the mainline is equal to assuming a fixed speed regarding to their impacts on the optimal result.

According to known researches and field observations, the fixed speed is seldom satisfied in real life. Variable Speed Limit could make vehicles' speeds relatively equally distributed to some extend but the cost would be the decrease of flow capacity in the mainstream.

If dynamic travel times on links are adopted, some more complex traffic flow models are required, e.g. CTM or AMOC model (Daganzo, 1994 and 1995; Kotsialos and Papageorgiou, 2004). This adoption may bring some benefits; but the undesirable impact on our analysis is obvious. All the models are only approximation to real life. According to Papageorgiou (1998), no model is perfect in its accuracy. Further surveillance of real world and adjustment of the model utilized are unavoidably required. This uncertainty and complexity of modeling and calibrating make it very hard, if possible, to answer our problems in a clear way. Even if we can use a more complicated model, it is doubtable that the result obtained in this way will be better than the one reached by adopting some simple and concise relationship assumptions. Theoretical analysis of a simple model may bring us some insights into the issues in question more easily.

Third, we assume that only one exit at the downstream end of the freeway is considered. No off-ramp will be considered in our model. Off-ramps may be considered in the future to improve the performance with respect to the holding capacities of partial mainline. To direct portion of flow to use the downstream off-ramps will decrease the pressure of queuing at the upstream on-ramps. This assumption also makes us avoid carrying out route choice. The reasonable explanation of this assumption lies in the fact that during an emergency evacuation, a downstream shelter always becomes the common destination for most of the evacuees using the freeway in question, so that the flow leaving from upstream off-ramps is insignificantly smaller. When some off-ramps need to be used to improve the flexibility of the model, a further study is required.

Fourth, we assume that the actual traffic demand is known or in a crude form of uncertainty without any related information about probability distribution. This means we will consider the deterministic demand with specific arrival patterns or uncertain demand in a crude form. An accurate dynamic demand does not exist according to Zhang and Levinson (2004). The assumption of demand linearly depending on time is not a proper assumption. Some studies take the assumption that the demand information is unavailable. But other limitations may show up in these studies, e.g. in So and Daganzo (2010) this may limit the final utilized ramps according to their physical features and in Zhang and Levinson (2004) an assumption of steady outflow at off-ramps is further required.

The fifth assumption used is that the freeway is under complete control. Here complete control of freeway means the entering flow from any on-ramp is controlled and the traffic flow in any part of the mainline should be kept under the flow capacity of that part. The flow capacity of any segment needs to be set using available information. Variable speed signs and other control facilities will be used to guarantee the mainline flow under capacity. When the contra-flow scheme is considered, the corresponding mainline capacity will change accordingly. By the way, if extra capacities of some stretches of mainline used under the free flow assumption are existed, we can identify them, and then make use of the extra-capacity to hold more vehicles so as to lighten the queuing burden on upstream on-ramps.

The last assumption is that the surface streets connecting with the freeway can be used to hold vehicles in view of the limited storage capacities of on-ramps during an emergency evacuation. Due to the emergency situation and our systematic optimization objective, the queue length in some on-ramps may be larger during an emergency evacuation than in an ordinary situation. But we should be aware that an in-depth research about the coordination of surface street traffic control and ramp metering is required in the future. A corridor of evacuation including freeway and arterial roadway could be a meaningful target of research for future. Note that this assumption of limited storage capacity of on-ramp is considered by many existing researches about freeway control.

*3.2. Main notations*

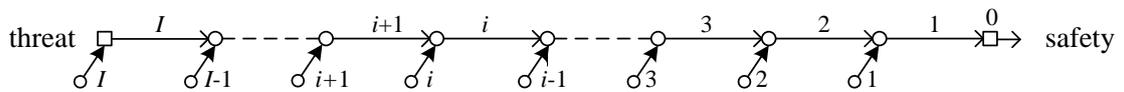

Fig. 1. A freeway used in an emergency evacuation.

We consider a freeway with the configuration as given above in Fig.1. Related main notations are as follows:

- $R$ is the set of on-ramps, indexed by $r \in R$;
- $L$ is the set of segments of mainline, indexed by $l \in L$;
- $\vec{T} = \{1, 2, 3, \cdots, K\}$ is the set of successive time intervals with equal span $T$, indexed by $k \in \vec{T}$ ($KT$ is the given evacuation time framework);
- $c_l$ is the maximum throughput in a typical time period $T$ for $l \in L$;

- $c_r$ is the maximum number of vehicles entering from $r \in R$ in a typical time period $T$;

- $f_r^k$ is the quantity of vehicles discharged from on-ramp $r$ during time interval $k$;

- $d_r^k$ is the number of vehicles arriving at on-ramp $r$ during time interval $k$;

- $D_r$ is the total demand at on-ramp $r$ such that $D_r = \sum_k d_r^k$ in the deterministic situation (when we consider an uncertain situation, $D_r$ becomes the upper bound of total demand during the feasible evacuation time at on-ramp $r$);

- $U_l$ is the set of on-ramps locating upstream of segment $l \in L$;

- $s_r$ is the storage (or holding) capacity of on-ramp $r \in R$.

*3.3. The basic LP models for a freeway evacuation*

There are two basic mathematical programming models which we will investigate below. The difference between them lies in the distinctive objective functions.

The first model aims at minimizing the total travel time of vehicles traversing the freeway. There are two ways to view our objective. One is to see it as all vehicles' travel times directly accumulated. When we assume that all the evacuees have arrived at corresponding on-ramps before the beginning of the evacuation and the free-flow speed is infinite, this travel time equals the total waiting time all vehicles experienced when to be driven away from freeway. The other way to view the following objective function is to see it as the summary of the throughput at the downstream exit weighted with the leaving time during the evacuation. Referring to Fig. 1, we can take our objective as to minimize the area bounded by the arrival curve, the departure curve and the horizontal line crossing point (0, D). Suppose that D is the total demand and time $t_D$ is the clearance time.

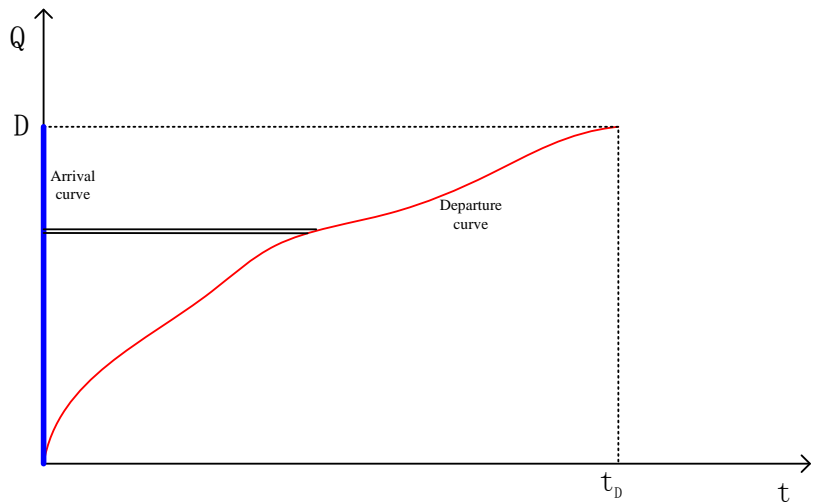

Fig. 2. Hypothetical arrival and departure curves for a controlled freeway.

The first specific mathematical programming model is given as follows:

$$[\text{MP1}]: \quad \min \sum_{r} \sum_{k} f_r^k kT \tag{1}$$

$$\sum_{r \in U_l} f_r^k \leq c_l, \quad \forall k \in \vec{T}, l \in L \tag{2}$$

$$f_r^k \leq c_r, \forall r, k \tag{3}$$

$$\sum_{k} f_r^k = D_r, \forall r \tag{4}$$

$$f_r^k \geq 0, \forall r, k \tag{5}$$

It is evident that above model is a linear programming model. Constraint (2) guarantees the mainline flow capacity cannot be exceeded. Constraint (3) makes sure the discharging flow from an on-ramp with an amount lower than the discharging capacity. Constraint (4) is used to guarantee all the evacuees can reach the downstream exit at the end of time horizon of evacuation. Constraint (5) is the nonnegative restriction on variables.

Note that in MP1 we only use $\sum_{k} f_r^k = D_r$ to connect discharging flow with the demand. This comes from the assumption that all the evacuees have arrived at their chosen on-ramps before the evacuation starts. So according to this assumption, the situation without evacuating vehicles at an on-ramp will not happen if there are still some vehicles which need to utilize this on-ramp in question. Later we will relax this assumption by adding new constraints.

Though to minimize the total travel time is important for an evacuation, the shortest time to finish the evacuation sometimes becomes the most desirable target, especially when the evacuation is time tight. Below a corresponding mathematical programming model is given to express above concern.

$$[\text{MP2}]: \quad \min zT$$

$$\min\{\tau \mid \sum_{k=1}^{\tau} f_r^k \geq D_r, \tau \in \vec{T}\} \leq z, \quad \forall r$$

$$\sum_{r \in U_l} f_r^k \leq c_l, \quad \forall k, l$$

$$f_r^k \leq c_r, \forall r, k \tag{6}$$

$$\sum_{k} f_r^k = D_r, \forall r$$

$$f_r^k \geq 0, \forall r, k$$

In MP2, the first constraint is a minimization sub-model used to determine the last time on an on-ramp when all the vehicles waiting in this ramp have entered the mainline. All the other constraints have the same meanings as in MP1. One important feature of MP2 is its nonlinear characteristic. Due to the large scale of the model, it is a real challenge to find its optimal solution sometimes. Though using a relative small freeway with 5 on-ramps and considering 100 time intervals in the time framework, we find out sometimes the commercial software Lingo11 fails to give out an optimal solution. By the way the computational time is generally unbearably long even lengthening to several hours to solve such a problem. Then a meaningful question, i.e. how to find out the optimal solution quickly, arises from above observation. Under some special conditions, e.g. suppose that any on-ramp has a sufficient discharging capacity to saturate its adjacent downstream mainline segment or the bottleneck only appears on the mainline, we can prove the set of optimal flow solutions to MP2 is equal to the one to MP1. The benefit resulted from this equivalent relation is the reduced computational time by using MP1 to replace MP2 to compute the minimum evacuation time.

But we should also notice that even without above equivalent relation, the minimum evacuation time for the whole system can be obtained by solving MP1. The way to do it is to appropriately adjust the evacuation time framework until $KT$ cannot be reduced further under the situation that a feasible solution can be obtained.

In a word, we should avoid solving a nonlinear programming model with a large scale. So let us just focus on MP1 from now on. Though it is easy to formulate a freeway system with fixed travel time on any segments of freeway mainline by the synchronization technique given by Zhang and Levinson (2004), the operation will make the subsequent analysis more complicated and the uncovering of relevant insights more difficult. In view of the uncertainty of real situations, it is unwise to assume fixed travel times for all the links sometimes.

*3.4. Extended model*

Based on the preceding basic models, we can extend them to deal with more complicated situations. First, we will consider embodying some predetermined priorities in the extended model. The priority may come from various considerations, e.g. HOV or bus priority on-ramps and risk levels changing down the freeway during an evacuation. A general objective can have the form of $\min \sum_r \sum_k w_r f_r^k kT$, where $w_r$ is the weight for on-ramp $r \in R$. For the MIN type of objective function, the smaller a weight is, the higher a priority will be given to the corresponding ramp.

Four types of evacuation strategies should be distinguished with respect to ramp discharging priority. The first is that any ramp has relative high flow-releasing priority over its downstream ramps. This can be called Innermost First Out (InFO) strategy (SO and Daganzo, 2010). The second one gives relative high priority to an on-ramp over its upstream ramps. After a careful observation, it is evident this type is equal to a freeway system without ramp metering during an emergency evacuation. Third, we can specify priorities to some on-ramps with special characters, e.g. with bus priority. The final type is an equity version that is $w_r = C, \forall r$ where $C$ is a given positive constant. Note that we do not consider the ramp metering explicitly here. It is obvious after ramp metering some evacuees may be delayed further comparing to the unmetered situation.

Note that we can further extend the weights to the time-variant form $w_r^k$. In this situation, the objective function will be $\min \sum_r \sum_k w_r^k f_r^k kT$. The reason for taking such a modification relies on the potential that we could adjust the priority of ramp to time dependent situations so as to optimize the performance of the whole system.

Another situation that can be dealt with by assigning proper priorities to some on-ramps is the appearance of unmetered on-ramps down the freeway. To make the model more practical in this situation, we had better set the weight indicating the priority of an unmetered on-ramp to a relative small value. Later we will show how to determine a proper weight for an unmetered on-ramp so as to make the resulted flow discharging pattern in accordance with reality.

Efficiency and equity are two aspects that can be considered here by putting different weights on different on-ramps. Later through sensitivity analysis, we can adjust the related weights so as to satisfy some requirements asked by stake holders. This adjustment should lead to an expected flow pattern on some specified on-ramps. But in an emergency evacuation, we should give more attention to special needs people and some other special situations, e.g. bus priority. The issue of equity is not the first concern here comparing to efficiency. So the general objective is to reduce the risk faced by people and to improve the efficiency of the whole system.

In the preceding subsection, we suppose all the evacuees have arrived at their chosen on-ramps before the implementation of an evacuation management strategy. But in many cases, the demand distribution is so different from that in the ordinary situation that no deterministic or stochastic demand estimate can be guaranteed in accordance with reality. One direct impact of above concern on our formulation is that the limitation of queuing at any on-ramp should be taken into account. A constraint such as $\sum_{\tau=1}^{k} f_r^\tau \leq \sum_{\tau=1}^{k} d_r^\tau$, $\forall r, k$ should be included in a practical model though it is hard to know the exact number $d_r^k$ of arrival vehicles to ramp $r$ during time interval $k$. In Section 5, we will consider a crude form of uncertain demand in more detail. The common utilized pattern of arrival rate during an emergency evacuation will also be presented in Section 5. An in-depth numerical analysis on the uncertain demand will be carried out in Section 6.

With above new concerns, we see that the objective $\min \sum_r \sum_k f_r^k kT$ of the basic model becomes improper now because this objective is only the total spent time for the throughput at the only exit under the assumption that every evacuee has reached his/her chosen on-ramp before the evacuation starts. Two equivalent objectives, i.e. $\min \sum_r \sum_k [\sum_{\tau=1}^{k} (d_r^\tau - f_r^\tau)T]$ and $\max \sum_r \sum_k [T(K+1-k)f_r^k]$, can be used now as the total waiting time under the situation with the information about time dependent $d_r^k$. The most impressive common characteristic of these two objective functions is their linearity. This property will benefit the subsequent analysis of a controlled freeway. Fig. 3 illustrates the equivalency between these two objectives. Referring to Fig. 3, we can relate the area bounded by the arrival curve, the departure curve and horizontal line crossing point (0, D) to the objective function $\min \sum_r \sum_k [\sum_{\tau=1}^{k} (d_r^\tau - f_r^\tau)T]$. The area bounded by departure curve, vertical line crossing point ($t_D$, 0) and the horizontal coordinate axis is corresponding to $\max \sum_r \sum_k [T(K+1-k)f_r^k]$. It is evident that with the fixed arrival curve the two objective functions are equivalent.

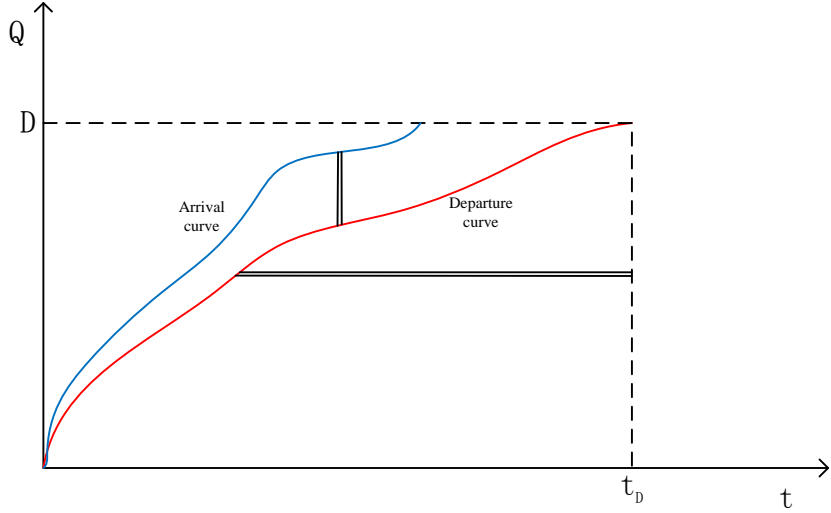

Fig. 3. Hypothetical arrival and departure curves for a controlled freeway with time dependent demand.

After considering the arrival pattern, we also need to take into account the storage capacity of on-ramp. To avoid the possible interference with surface traffic due to the long queue w.r.t. an on-ramp, a constraint such as $\sum_{\tau=1}^{k} d_r^{\tau} - \sum_{\tau=1}^{k} f_r^{\tau} \leq s_r, \forall r \in R, k \in \vec{T}$ should be added into the extended model. The constraint in conjunction with the constraint $\sum_{\tau=1}^{k} f_r^{\tau} \leq \sum_{\tau=1}^{k} d_r^{\tau} \quad \forall r, k$ can be viewed as a simplified version of METANET related queue constraints (Kotsialos and Papageorgiou, 2004). The simplification of above constraints comparing with METANET comes from the assumption of a free flow state of mainstream. In METANET, the ramp discharging flow is limited by the current mainstream state.

If the first responders who may need to use the freeway to reach the site of risk at the beginning of an evacuation are considered, some control methods, e.g. appointing some contra-flow lanes, will be adopted to facilitate their rescue action. But with time, the lanes used by first responders will gradually be reused by ordinary evacuees. So the capacity of part of the freeway is time variant. In view of this observation, the constraint of $\sum_{r \in U_l} f_r^k \leq c_l$ should be changed to be $\sum_{r \in U_l} f_r^k \leq c_l^k$. How to make use of the dynamic changing extra mainline capacities on different stretches is not only a challenge but also an interesting research issue. The same logic of time variant parameters can be applied to storage capacity $s_r$, too.

With all above mentioned modification, the final extended model is given as follows:

[MP3]:  $\quad \max \sum_r \sum_k [(K+1-k) w_r^k f_r^k] \quad$ (or $\quad \min \sum_r \sum_k [\sum_{\tau=1}^{k} w_r^{\tau} (d_r^{\tau} - f_r^{\tau}) T])$ (7)

S.T. $\quad \sum_{r \in U_l} f_r^k \leq c_l^k, \quad \forall k \in \vec{T}, l \in L$ (8)

$\quad f_r^k \leq c_r, \forall r \in R, k \in \vec{T}$ (9)

$\quad \sum_k f_r^k \geq D_r, \forall r \in R$ (10)

$\quad f_r^k \geq 0, \forall r \in R, k \in \vec{T}$ (11)

$$\sum_{\tau=1}^{k} f_r^{\tau} \leq \sum_{\tau=1}^{k} d_r^{\tau}, \forall r \in R, k \in \vec{T} \tag{12}$$

$$\sum_{\tau=1}^{k} d_r^{\tau} - \sum_{\tau=1}^{k} f_r^{\tau} \leq s_r^k, \forall r \in R, k \in \vec{T} \tag{13}$$

Note that we change the equation $\sum_k f_r^k = D_r$ to inequality $\sum_r f_r^k \geq D_r$ in the extended model. In view of the relation $D_r = \sum_k d_r^k$ and constraint (12), this change will not influence the final solution of MP3. This modification will make the subsequent analysis of uncertain demand using AARC become easier.

Before we go on to the sensitivity analysis of above model, another interesting question arising from above model is worthy of further consideration. That is how to judge whether the MP3 has a feasible solution. It is very important for evacuation managers. A positive answer to this question is a prerequisite for evacuation agents to design an effective control strategy for managing freeway evacuation. The constraints including the limited holding capacity on a ramp, the time variant arrival flow pattern to on-ramps, the changing mainline traffic capacity, and the discharging capacities related to on-ramps, make it difficult to judge whether a feasible solution exists to the extended model. Fortunately, the linearity of the model in conjunction with powerful commerce software existing saves us from this difficult situation. We can know the answer simply through solving it with any suitable software on hand.

## 4. The sensitivity analysis

In view of the existence of many powerful solution methods nowadays, e.g. column generation algorithm, we will just focus on the application of related results about sensitivity analysis of LP. We also suppose readers will choose some commercial software to carry out their analysis, e.g. LINGO11. For simplification of expression, we will take LINGO11 as our working software. Substantial useful information, e.g. dual price accessible after the solving of our model, can be used to facilitate our analysis.

In this section, we focus on physical explanations and possible applications of various sensitivity analyses. The related theoretical proofs about sensitivity analysis can be found in many proper textbook, e.g. Wayne L. Winston (2004). The specific analytical examples will be given in Section 5 later. With respect to sensitivity analysis used to prepare and respond an evacuation, there are three very important aspects we will investigate below.

First, let us consider how to realize the predetermined priorities assigned to some specified on-ramps. In view of that our model is to minimize the total waiting time, to embody the predetermined priorities in the objective function, we have added the weights $w_r^k$ in the objective function. To consider the minimization type of objective function, it is evident the smaller the weight is, the higher the priority will be. To make the predetermined priority possible to realize, for an ordinary on-ramp a relative high value of its corresponding weight should be set to. Given a series of weights, after the solving, the predetermined priorities may be realized or not. When the later happens, we should adjust the weights in a proper way. For example we observe one ramp with higher priority releases no flow during a time period. We expect this on-ramp can release more flow during this period with the slight impact on the performance of the whole system. An intuitive way to realize it is to reduce the corresponding weights (here we consider the minimization type of objective function) for the on-ramp in this time period, and then to resolve the model again to observe if our aim is achieved. During this process, the initial results supplied by LINGO11 about the allowable changing ranges of the corresponding objective function coefficients will be helpful to give us the information to adjust the weights. The current values of related variables corresponding to these weights also give us the general information about the change of objective function value. Practitioners can estimate the lost after this adjustment of weights.

Through this process, we may find out the possible correlated on-ramps that can discharge flow interchangeably without changing the clearance time of evacuation network. It means to let one on-ramp release flow, the other on-ramp or ramps in the group of these correlated ramps will release no flow during some specific time period. One important issue arising here is whether there are a series of optimal weights which by adopting we can realize our predetermined priorities and attain the optimal performance of the whole freeway system at the same time. It is not easy to answer this question at present. To give a satisfying answer, an optimal control problem needs to be formulated and resolved. We think an in-depth investigation about this concern is required in the near future.

The second issue we are concerned with is how to use the sensitivity analysis of LP to make use of the mainline traffic capacity. In our model, the freeway mainline has been divided into many segments. Any one of these segments has its uniform geological and physical features. As mentioned earlier, the capacity $c_l$ of segment $l$ may be changeable during an emergency evacuation due to specific usage in some period or possible incidents. There are two basic situations about the use of mainline segment capacity. On the one hand, some segments may be under its capacity used during an emergency evacuation; on the other hand, some segments' capacities may be used up. For the former, we expect to take advantage of the extra capacities to hold more vehicles in the mainline to lighten the pressure of queuing in upstream on-ramps. For the latter, the segments with used up capacities, as the critical dynamic bottlenecks to an efficient emergency evacuation, should be given more attentions to for keeping them functioning smoothly and safely. The related information can be found in the report of LINGO11 in the data of Slack or Surplus corresponding to the constraints of mainline capacities. A positive value of slack or surplus means the existence of extra capacity in the corresponding segment. A slack or surplus of zero tells us the capacity of the segment in question has been used up. The allowable changing ranges of RHS related to mainline segments can be used to analyze the impact of an adjusted capacity on the objective value. Like the sensitivity analysis of weights, we can see the influence of adjusted RHS by solving the LP model repeatedly.

A valuable conclusion given in existing literature tells us that a relative steady occupancy (or density) in a segment is very likely existing comparing to a uniform average speed of vehicles. If we can improve the occupancy and at the same time keep the mainstream running smoothly, we obtain an ideal state of freeway system to implement an efficient evacuation and to hold more vehicles downstream with relative low risk level.

If we know the fundamental relationship between traffic flow and density on an unsaturated segment $l$ as illustrated in Fig. 4, we can roughly estimate the possible extra holding capacity in this segment. Suppose that the capacity $c_l$ is equal to $q_{critical}$ and current flow corresponds to point A in the figure. We can expect to improve the density from $k_A$ to $k_B$ by some control method, e.g. variable speed limits, so as to buffer more vehicles in this segment of mainline and at the same time keep an adequate throughput. The rough quantity of vehicles to be buffered can be calculated by $(k_B - k_A) \text{L}_l \lambda_l$ where $\text{L}_l$ and $\lambda_l$ are the length and the number of lanes of segment $l$, respectively. These preliminary analytical results can help the subsequent simulation experiment with more complex traffic flow model as basis, e.g. CTM. As is well known, general LP models based on CTM may hold flows through improving the maximum capacity of a cell and lowering the speed of related vehicles. Based on the result of our model, we can realize the expected buffering in part of mainline through adjusting the maximum capacities of cells composing of the related segments. This promising research is under way now. We expect to support an in-depth report later.

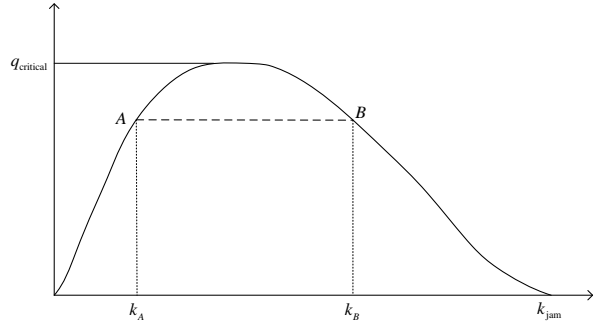

Fig. 4. The fundamental relationship between flow $q$ and density $k$.

The third issue we want to address is about how to set a proper value to the storage capacity $s_r^k$ of an on-ramp $r$. For an actual evacuation using a freeway, the possible interference with surface road traffic due to inadequate storage capacity of on-ramp must be considered carefully. Two reasons lies under this interference. One is the abruptly increased uncertain demand during a short time to an on-ramp. The other is the improper ramp-metering rates. Given the demand pattern and a series of discharging weights, through solving MP3 with increased $s_r^k$, we can check whether some constraints with respect to on-ramp storage capacities have been violated. If it is true, we have two ways to deal with this situation. One is to adjust the priorities of ramp discharging. This way may make use of some unsaturated ramps to buffer more vehicles, and at the same time make the on-ramp in question to release flow as quickly as possible. The other choice to solve this problem is to consider the possible available space of adjacent surface roads to buffer the incoming flow to the ramp in question. New increased storage capacity can be coded in MP3 to find the new optimal solution corresponding to the new parameters. Comparing with the traditional methods using complex traffic flow model and simulation technique, our method is easy and fast. In view of the emergency situation, it is a reliable way to fix the problem of shortage of storage space.

## 5. The robust optimization model to deal with demand uncertainty

Different from sensitivity analysis, robust optimization is to search the optimal solution to the model with uncertainty constraints in a crude form. This kind of analysis is especially desired in an emergency evacuation situation. As is well known, it is very difficult to estimate or forecast an accurate traffic demand for an emergency evacuation. So in such an uncertain situation, to assume some crude form of demand is not only a reasonable choice but also a possible practice most of time resulted from readily available information.

*5.1. The uncertainty of dynamic demand*

In MP3, the demand information as RHS of constraints is to influence the result of the model. We suppose that the exact amount of demand $d_r^k$ is difficult to know before time $k$, but the possible changing range and the upper bound of the accumulative demand during the whole time horizon are known. The specific polyhedral constraints are given as follows:

$$\sum_{k\in \bar{T}} d_r^k \le D_r \text{ and } \underline{d}_r^k \le d_r^k \le \bar{d}_r^k, \forall r,k. \tag{14}$$

To make the subsequent expressions concise, we denote a set to embody the above constraints, that is

$$d_r^k \in U_r^k \equiv \{d_r^k : \underline{d}_r^k \leq d_r^k \leq \overline{d}_r^k, \sum_{k \in \vec{T}} d_r^k \leq D_r\}. \tag{15}$$

Except above polyhedral form of uncertainty set, another commonly used uncertainty set has a simple box form as follows:

$$d_r^k \in U_r^k \equiv \{d_r^k : \tilde{d}_r^k(1-\theta) \leq d_r^k \leq \tilde{d}_r^k(1+\theta)\}, \tag{16}$$

where $\theta$ denotes uncertainty level and $\tilde{d}_r^k$ is the nominal demand in on-ramp $r$ during time interval $k$. It is easy to embed (16) in (15).

With above crude uncertainty constraints, MP3 will become an intractable semi-infinite dimensional mathematical programming. Without further treatment, it is hard to obtain its optimal solution, sometime even a feasible solution.

Though the pattern of dynamic demands during an emergency evacuation is often scenario-dependent, one famous deterministic arrival pattern is broadly accepted in practical and theoretical field. Jamei B. (1984) started to use mobilization curve, i.e. S-shaped curve of accumulative arrival rate, to describe the arrival pattern of evacuees at an appointed location. For the convenience of the subsequent numerical analysis in next section, we will give a brief introduction of this pattern. The mobilization curve is represented by the following equation:

$$\gamma_r^k = 1/\{1+\exp[-LR(k-HF)]\} \tag{17}$$

where $\gamma_r^k$ is the cumulative percentage of evacuees arrived at on-ramp $r$ in a freeway network by time $k$, $LR$ is the response rate of public to the disaster order, and $HF$ is the half loading time. Generally speaking, the loading rate $LR$ is dependent on the efficiency of communicating the evacuation order to citizens, the response level of citizens to the evacuation order, and the severity of the disaster (Jamei B., 1984). According to (17), $HF$ is the time instant when the cumulative percentage reaches 50%. A large value of **HF** means that many evacuees may choose to delay their response to an evacuation order. Fig. 4 will give reader some intuitive impression about this deterministic arrival pattern.

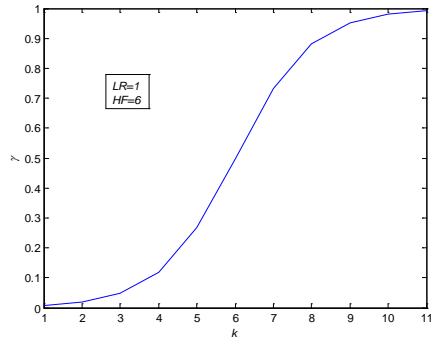

Fig.5. Mobilization curve with *LR*=1 and *HF*=6.

*5.2. AARC model*

To construct a tractable model, we suppose that an affine relation between the discharging flow $f_r^k$ and the demands $d_s^\tau, \forall s \in R, \tau \in T_k$, exists. Set $T_k = \{1,\cdots,k\}$ is a set of time intervals. The following specific relation

exists:

$$f_r^k = \eta_{rk}^{-1} + \sum_{s\in R}\sum_{\tau\in T_k} \eta_{rk}^{s\tau} d_s^\tau. \tag{18}$$

In Eq. (18), $\eta_{rk}^{-1}$ and $\eta_{rk}^{s\tau}$ are non-adjustable variables. Using the expression in Eq. (18) to substitute $f_r^k$ in MP3, we can transform MP3 into the following form:

$$\begin{aligned}
&\min z \quad \text{(M-AARC)} \\
&\sum_{r\in R}\sum_{k\in \vec{T}}[(k-1-K)w_r^k(\eta_{rk}^{-1}+\sum_{s\in R}\sum_{\tau\in T_k}\eta_{rk}^{s\tau}d_s^\tau)]\le z, \forall d_r^k \in U_r^k \\
&\sum_{r\in U_l}(\eta_{rk}^{-1}+\sum_{s\in R}\sum_{\tau\in T_k}\eta_{rk}^{s\tau}d_s^\tau)\le c_l^k, \ \ \forall d_r^k \in U_r^k, k\in \vec{T}, l\in L \\
&\eta_{rk}^{-1}+\sum_{s\in R}\sum_{\tau\in T_k}\eta_{rk}^{s\tau}d_s^\tau \le c_r, \forall d_r^k \in U_r^k, k\in \vec{T}, r\in R \\
&\sum_{k\in \vec{T}} d_r^k - \sum_{k\in \vec{T}}(\eta_{rk}^{-1}+\sum_{s\in R}\sum_{\tau\in T_k}\eta_{rk}^{s\tau}d_s^\tau)\le 0, \forall d_r^k \in U_r^k, r\in R \\
&\sum_{v\in T_k}(\eta_{rv}^{-1}+\sum_{s\in R}\sum_{\tau\in T_v}\eta_{rv}^{s\tau}d_s^\tau)-\sum_{v\in T_k}d_r^v \le 0, \ \ \forall d_r^k \in U_r^k, k\in \vec{T}, r\in R \\
&\sum_{v\in T_k}d_r^v - \sum_{v\in T_k}(\eta_{rv}^{-1}+\sum_{s\in R}\sum_{\tau\in T_v}\eta_{rv}^{s\tau}d_s^\tau)\le s_r^k, \forall d_r^k \in U_r^k, k\in \vec{T}, r\in R \\
&\eta_{rk}^{-1}+\sum_{s\in R}\sum_{\tau\in T_k}\eta_{rk}^{s\tau}d_s^\tau \ge 0, \forall d_r^k \in U_r^k, k\in \vec{T}, r\in R
\end{aligned} \tag{19}$$

After above substitution, the new linear programming model has $\eta_{rk}^{-1}$ and $\eta_{rk}^{s\tau}$ as its decision variables. The original variables $f$ can be calculated using the linear relation of Eq. (18). By solving M-AARC (19), the objective value denoted as $z_{M-ARRC}^*$ is the guaranteed upper bound value for all realization of uncertain data under the assumption of linear dependency. This value can also be viewed as the optimal estimate of the weighed total waiting time in the worst case.

M-AARC (19) remains semi-infinite and intractable. Some further modification is still required. This can be done by using the strong dual theorem of LP to change the constraints. The final tractable formulation is presented below.

$$\begin{aligned}
&\min z \quad \text{(M-AARC1)} \\
&\sum_{\tau\in \vec{T}}\sum_{s\in R}(\bar{d}_s^\tau \lambda_{\tau s}^{11} - \underline{d}_s^\tau \lambda_{\tau s}^{12}) + \sum_{s\in R} D_s \lambda_s^{13} \le z - \sum_{r\in R}\sum_{k\in \vec{T}}(k-1-K)w_r^k \eta_{rk}^{-1} \\
&\lambda_{\tau s}^{11}-\lambda_{\tau s}^{12}+\lambda_s^{13} = \sum_{k=\{\tau\ldots K\}}\sum_{r\in R}(k-1-K)w_r^k \eta_{rk}^{s\tau}, \ \ \forall \tau \in \vec{T}, s\in R \\
&\lambda_{\tau s}^{11},\lambda_{\tau s}^{12},\lambda_s^{13} \ge 0, \ \ \forall \tau \in \vec{T}, s\in R \\
&\sum_{\tau\in \vec{T}}\sum_{s\in R}(\bar{d}_s^\tau \lambda_{lk\tau s}^{21}-\underline{d}_s^\tau \lambda_{lk\tau s}^{22})+\sum_{s\in R}D_s \lambda_{lks}^{23} \le c_l^k - \sum_{r\in U_l}\eta_{rk}^{-1}, \forall k\in \vec{T}, l\in L
\end{aligned}$$

$$\lambda_{lk\tau s}^{21} - \lambda_{lk\tau s}^{22} + \lambda_{lks}^{23} = \sum_{r \in U_l} \eta_{rk}^{s\tau}, \quad \forall \tau = \{1\ldots k\}, s \in R, k \in \vec{T}, l \in L$$

$$\lambda_{lk\tau s}^{21} - \lambda_{lk\tau s}^{22} + \lambda_{lks}^{23} = 0, \quad \forall \tau = \{k+1\ldots K\}, s \in R, k \in \vec{T}, l \in L$$

$$\lambda_{lk\tau s}^{21}, \lambda_{lk\tau s}^{22}, \lambda_{lks}^{23} \geq 0, \quad \forall \tau \in \vec{T}, s \in R, k \in \vec{T}, l \in L$$

$$\sum_{\tau \in \vec{T}} \sum_{s \in R} (\overline{d}_s^\tau \lambda_{rk\tau s}^{31} - \underline{d}_s^\tau \lambda_{rk\tau s}^{32}) + \sum_{s \in R} D_s \lambda_{rks}^{33} \leq c_r - \eta_{rk}^{-1}, \forall k \in \vec{T}, r \in R$$

$$\lambda_{rk\tau s}^{31} - \lambda_{rk\tau s}^{32} + \lambda_{rks}^{33} = \eta_{rk}^{s\tau}, \quad \forall \tau = \{1\ldots k\}, s \in R, k \in \vec{T}, r \in R$$

$$\lambda_{rk\tau s}^{31} - \lambda_{rk\tau s}^{32} + \lambda_{rks}^{33} = 0, \quad \forall \tau = \{k+1\ldots K\}, s \in R, k \in \vec{T}, r \in R$$

$$\lambda_{rk\tau s}^{31}, \lambda_{rk\tau s}^{32}, \lambda_{rks}^{33} \geq 0, \quad \forall \tau \in \vec{T}, s \in R, k \in \vec{T}, r \in R$$

$$\sum_{\tau \in \vec{T}} \sum_{s \in R} (\overline{d}_s^\tau \lambda_{r\tau s}^{41} - \underline{d}_s^\tau \lambda_{r\tau s}^{42}) + \sum_{s \in R} D_s \lambda_{rs}^{43} \leq \sum_{k \in \vec{T}} \eta_{rk}^{-1}, \forall r \in R$$

$$\lambda_{r\tau s}^{41} - \lambda_{r\tau s}^{42} + \lambda_{rs}^{43} = I_{(r=s)} - \sum_{k \in \vec{T}} \eta_{rk}^{s\tau} I_{(\tau \leq k)}, \quad \forall \tau \in \vec{T}, s \in R, r \in R$$

$$\lambda_{r\tau s}^{41}, \lambda_{r\tau s}^{42}, \lambda_{rs}^{43} \geq 0, \quad \forall \tau \in \vec{T}, s \in R, r \in R$$

$$\sum_{\tau \in \vec{T}} \sum_{s \in R} (\overline{d}_s^\tau \lambda_{rk\tau s}^{51} - \underline{d}_s^\tau \lambda_{rk\tau s}^{52}) + \sum_{s \in R} D_s \lambda_{rks}^{53} \leq -\sum_{\upsilon=1}^{k} \eta_{r\upsilon}^{-1}, \forall k \in \vec{T}, r \in R$$

$$\lambda_{rk\tau s}^{51} - \lambda_{rk\tau s}^{52} + \lambda_{rks}^{53} = \sum_{\upsilon=1}^{k} \eta_{r\upsilon}^{s\tau} I_{(\tau \leq \upsilon)} - I_{(r=s, \tau \leq k)}, \quad \forall \tau \in \vec{T}, s \in R, k \in \vec{T}, r \in R$$

$$\lambda_{rk\tau s}^{51}, \lambda_{rk\tau s}^{52}, \lambda_{rks}^{53} \geq 0, \quad \forall \tau \in \vec{T}, s \in R, k \in \vec{T}, r \in R$$

$$\sum_{\tau \in \vec{T}} \sum_{s \in R} (\overline{d}_s^\tau \lambda_{rk\tau s}^{61} - \underline{d}_s^\tau \lambda_{rk\tau s}^{62}) + \sum_{s \in R} D_s \lambda_{rks}^{63} \leq s_r^k + \sum_{\upsilon=1}^{k} \eta_{r\upsilon}^{-1}, \forall k \in \vec{T}, r \in R$$

$$\lambda_{rk\tau s}^{61} - \lambda_{rk\tau s}^{62} + \lambda_{rks}^{63} = I_{(r=s, \tau \leq k)} - \sum_{\upsilon=1}^{k} \eta_{r\upsilon}^{s\tau} I_{(\tau \leq \upsilon)}, \quad \forall \tau \in \vec{T}, s \in R, k \in \vec{T}, r \in R$$

$$\lambda_{rk\tau s}^{61}, \lambda_{rk\tau s}^{62}, \lambda_{rks}^{63} \geq 0, \quad \forall \tau \in \vec{T}, s \in R, k \in \vec{T}, r \in R$$

$$\sum_{\tau \in \vec{T}} \sum_{s \in R} (\overline{d}_s^\tau \lambda_{rk\tau s}^{71} - \underline{d}_s^\tau \lambda_{rk\tau s}^{72}) + \sum_{s \in R} D_s \lambda_{rks}^{73} \leq \eta_{rk}^{-1}, \forall k \in \vec{T}, r \in R$$

$$\lambda_{rk\tau s}^{71} - \lambda_{rk\tau s}^{72} + \lambda_{rks}^{73} = -\eta_{rk}^{s\tau}, \quad \forall \tau = \{1\ldots k\}, s \in R, k \in \vec{T}, r \in R$$

$$\lambda_{rk\tau s}^{71} - \lambda_{rk\tau s}^{72} + \lambda_{rks}^{73} = 0, \quad \forall \tau = \{k+1\ldots K\}, s \in R, k \in \vec{T}, r \in R$$

$$\lambda_{rk\tau s}^{71}, \lambda_{rk\tau s}^{72}, \lambda_{rks}^{73} \geq 0, \quad \forall \tau \in \vec{T}, s \in R, k \in \vec{T}, r \in R \tag{20}$$

Note that $\lambda$ is a set of dual variables and the numerical indexes are used for notational simplicity. In above model, $I_{(s)}$ is defined as a binary number. When the condition $s$ is satisfied, $I_{(s)}$ will take value 1; otherwise, it will take value 0.

**Theorem 1**: Given polyhedral uncertainty set, $U_r^k$, the affinely adjustable robust counterpart of MP3 becomes the LP problem (20) and thus computationally tractable.

**Proof**. By using the following relationship, we can reformulate each constraint affected by uncertain data as an equivalent LP problem,

$$\sum_{\tau \in T_k} \alpha_r^\tau d_r^\tau \leq v \quad \forall d_r^k \in U_r^k = \{\underline{d}_r^k \leq d_r^k \leq \overline{d}_r^k, \sum_{k \in \vec{T}} d_r^k \leq D_r\} \Leftrightarrow \max_{d_r^\tau}(\sum_{\tau \in T_k} \alpha_r^\tau d_r^\tau) \leq v.$$

Without loss of generality, the polyhedral uncertainty set is represented as $Ad \leq b$ and $\max_{d_r^\tau}(\sum_{\tau \in T_k} \alpha_r^\tau d_r^\tau)$ is written as (P). By applying strong duality property, we can derive an equivalent constraint with dual problem (D) on the left hand side.

$$\begin{array}{l} \max \alpha d \leq v \\ s.t. \\ Ad \leq b \\ d \text{ urs} \end{array} \Leftrightarrow \begin{array}{l} \min b\lambda \leq v \\ s.t. \\ A^T \lambda = \alpha \\ \lambda \geq 0 \end{array}.$$

Therefore, there exists $\lambda$ satisfying $b\lambda \leq v$, $A^T\lambda = \alpha$ and $\lambda \geq 0$ and the equivalent AARC of MP3 becomes tractable. □

In practice, we can solve the AARC model (20) to obtain a series of $\eta_{rk}^{-1}$ and $\eta_{rk}^{s\tau}$ off-line, and then use (18) to get the final discharging flow strategy on-line. Another choice in practice is to roll forward the time horizon and solve the model repeatedly with the already revealed demand information. In view of the linearity of model (20), the latter will not bring us more computation complexity.

## 6. Numerical analysis

In this section, we will use simple examples to illustrate the effectiveness and efficiency of main ideas given in Section 4 and 5. We consider a freeway with 5 on-ramps and 5 mainline segments which has the similar structure as shown in Fig. 1.

*6.1. Basic parameters for our numerical experiment*

Two types of deterministic demands will be investigated in subsection 6.2. One follows the S-shaped curve. The other simply assumes a fixed arrival rate for each of on-ramps. For both types, the time horizon is set to 20. The effective arrival time limits are set to 10 and 14 with respect to fixed arrival rates and S-shaped curve, respectively. Time horizon means the evacuation would be finished in the specified time periods, and the evacuees will arrive at their chosen on-ramps before the specified effective arrival time limits.

For the ordinary benchmark situation, we set all weights to 1 since we only consider the MAX type of objective function from now on. The information of various capacities for the benchmark scenario is summarized in Table 1.

Note that we assume a series of time-invariant capacities for this scenario.

**Table 1**

The basic parameters about this freeway in the ordinary situation.

| No. of camps (or links) | 1 | 2 | 3 | 4 | 5 |
|---|---|---|---|---|---|
| $c_r$ | 4 | 5 | 3 | 4 | 3 |
| $c_l$ | 14 | 14 | 16 | 10 | 10 |
| $s_r$ | 9 | 8 | 10 | 8 | 6 |
| $d_r$ | 4 | 3 | 4 | 4 | 3 |

*6.2. Sensitivity analysis to address priorities and capacities*

*6.2.1. Fixed arrival rates*

Through solving the MP3, we get the discharging flow pattern for the ordinary scenario presented in Table 2. The objective value is 2531 and the clearance time is 14. There are total 120 variables and 426 constraints in this model. The corresponding cumulative exiting flow is shown in Fig. 6. From the figure, we can see the system stays on its bottleneck capacity from the beginning to time interval 12. Only in the final two time intervals near the end of clearance the exiting flow decreases.

**Table 2**

The discharging flows in the clearance time for the ordinary situation.

| Time interval | 1 | 2 | 3 | 4 | 5 | 6 | 7 | 8 | 9 | 10 | 11 | 12 | 13 | 14 |
|---|---|---|---|---|---|---|---|---|---|---|---|---|---|---|
| Ramp 1 | 1 | 4 | 3 | 2 | 4 | 4 | 1 | 4 | 4 | 4 | 4 | 1 | 4 | 0 |
| Ramp 2 | 3 | 0 | 2 | 5 | 5 | 0 | 5 | 2 | 0 | 0 | 0 | 4 | 4 | 0 |
| Ramp 3 | 3 | 3 | 3 | 3 | 3 | 3 | 3 | 3 | 3 | 3 | 3 | 3 | 3 | 1 |
| Ramp 4 | 4 | 4 | 4 | 4 | 0 | 4 | 2 | 2 | 4 | 4 | 4 | 4 | 0 | 0 |
| Ramp 5 | 3 | 3 | 2 | 0 | 2 | 3 | 3 | 3 | 3 | 3 | 3 | 2 | 0 | 0 |

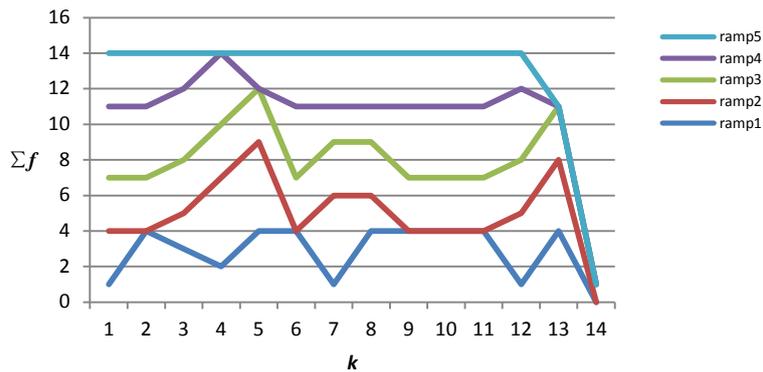

Fig. 6. The cumulative exiting flow in the ordinary situation.

If on-ramp 2 has bus exclusive lane, a higher priority may be assigned to it. In the ordinary situation, we find out that the flow discharging pattern is not satisfied for us. We expect more evacuees can leave earlier from on-ramp 2. To realize above expectation, we try to change the weight to on-ramp 2 to influence the final flow pattern. So we increase the weight to ramp 2 from 1 to 10 and resolve the modified model, the result is summarized in Fig. 7 and Table 3. The objective value increases to 6590 due to the renewed weights. The clearance time does not change after this modification. In Fig. 8, the discharging flow patterns to on-ramp 2 after and before above modification are compared. It is evident from Fig. 8 that the improved weight makes more evacuees can leave earlier from on-ramp 2. The comparison also shows us the possibility of changing specific discharging flow pattern without changing the clearance time.

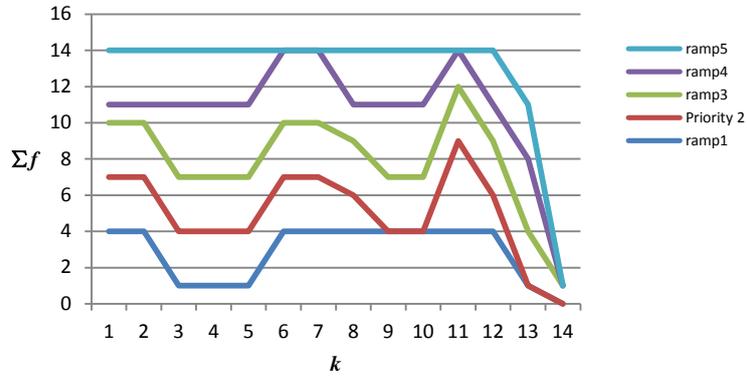

Fig.7. The cumulative flow curve with $w_2=10$.

**Table 3**

The discharging flows in the clearance time with $w_2=10$.

| Time interval | 1 | 2 | 3 | 4 | 5 | 6 | 7 | 8 | 9 | 10 | 11 | 12 | 13 | 14 |
|---|---|---|---|---|---|---|---|---|---|---|---|---|---|---|
| Ramp 1 | 4 | 4 | 1 | 1 | 1 | 4 | 4 | 4 | 4 | 4 | 4 | 4 | 1 | 0 |
| Ramp 2 | 3 | 3 | 3 | 3 | 3 | 3 | 3 | 2 | 0 | 0 | 5 | 2 | 0 | 0 |
| Ramp 3 | 3 | 3 | 3 | 3 | 3 | 3 | 3 | 3 | 3 | 3 | 3 | 3 | 3 | 1 |
| Ramp 4 | 1 | 1 | 4 | 4 | 4 | 4 | 4 | 2 | 4 | 4 | 2 | 2 | 4 | 0 |
| Ramp 5 | 3 | 3 | 3 | 3 | 3 | 0 | 0 | 3 | 3 | 3 | 0 | 3 | 3 | 0 |

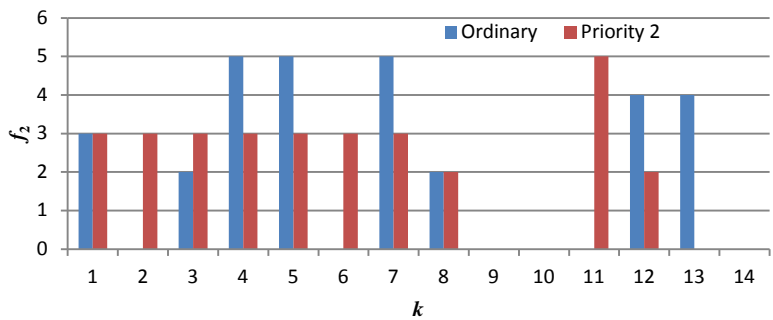

Fig.8. The comparison of discharging flow patterns of on-ramp 2.

To see if the changing weights can realize the InFO strategies, we set a weight $r^3$ to on-ramp $r$. Through solving the corresponding MP3, the cumulative released flow pattern is shown in Fig. 9. It is obvious that the higher weight has lead to earlier releasing of more evacuees from corresponding on-ramp. A similar analysis can be carried out to the downstream first scenario in which we set a weight $(9-r)^3$ to on-ramp $r$. The result is demonstrated by Fig. 10. Through comparing the figures 6, 7, 9 and 10, we can see that the weight-changing can attain the predetermined priorities in a convenient way at the same time without stretching the clearance time of system.

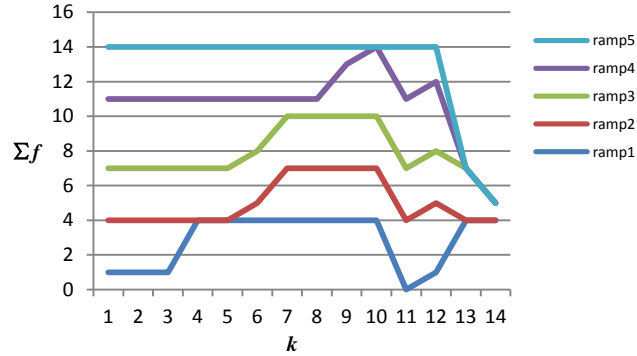

Fig.9. The cumulative flow curve for the Inner Most First Out scenario

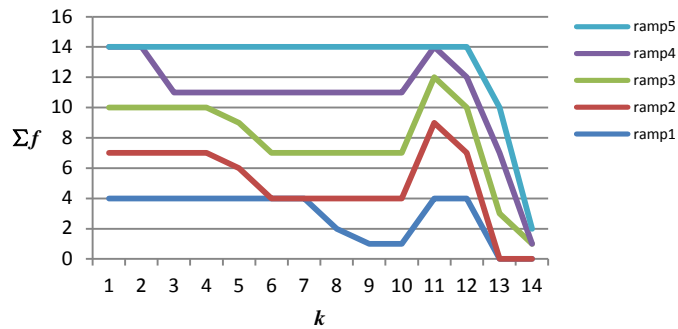

Fig.10. The cumulative flow curve for the Outer Most First Out scenario

After the analysis of weight-changing influence, let us consider the changing mainline capacity. Based on the information supplied by Lingo 11 after solving the model MP3 w.r.t. the ordinary situation, the extra capacities are 0, 1, 6, 3, and 7 on mainline segments 1, 2, 3, 4, and 5, respectively. So we can judge the bottleneck lies in mainline segment 1. Since the others with positive extra capacity during the whole evacuation, they could all be used as buffer zone to holding more vehicles. In Table 4, we list the impacts of changing bottleneck capacity on objective value and clearance time. We can see a little bit of decrease of bottleneck capacity will lead to non feasible solution to the model and infinite clearance time. The increase of bottleneck capacity can improve the objective value, but has no impact on the clearance time. The results show us that more evacuees will leave the freeway and reach the shelter earlier when the bottleneck or critical segment's capacity is increased. To make the whole system operate healthily, more care should be taken of the bottleneck of mainline.

**Table 4**

The impact of changing mainline capacity $\triangle c_1$ on objective value.

| $\triangle c_1$ | -0.2 | -0.1 | 0.1 | 0.5 | 1 |
|---|---|---|---|---|---|
| Objective Value | Non feasible solution | 2523.2 | 2538.8 | 2570 | 2605 |
| Clearance Time | Infinite | 14 | 14 | 14 | 14 |

During the evacuation, segment 3 keeps a high extra dynamic capacity, so we can make use of it to buffer more vehicles so as to lighten the queuing burden in on-ramps. Through calculating the modified models, we can see that the capacity of mainline segment 3 can be decreased to 8.6 with the same objective value (note that the specific discharging flows have changed). When the capacity of segment 3 is further decreased to 8.5, the model will produce non feasible solution. Note that when the changing is restricted by the allowable range given by Lingo, i.e. the allowable increase is infinite and the allowable decrease is 6, the objective value and the optimal discharging flows will both keep the same.

For the ordinary scenario, we find out that except on-ramp 5 with extra storage capacity 1, the others all used up their storage capacity at some time intervals during the evacuation. Further observation of the result tells us, for ramps 1 to 4 at time interval 10, they all used up their storage capacities. For any one of them only during a few time intervals its storage capacity was used up, e.g. for on-ramp 3 only in time interval 10 its capacity is used up. This observation suggests us a possible way to avoid the broken down of system during the critical time intervals that is temporary using the adjacent surface streets to hold more arrived vehicles. Suppose that we decrease the storage capacity of any one of ramps 1 to 4 during time intervals 9, 10 and 11 by 0.1, then let us see what will happen to the whole system. Surprisingly, we find out under any one of situations generated by decreasing one ramp storage capacity during time intervals from 9 to 11, the model has no feasible solution existing. So keeping an eye on the ramps during these time periods is very important. The possible cooperation between surface roads and on-ramps is required during such critical time intervals.

The other observation is that even if we increase all ramps' storage capacities by 2 during time intervals from 9 to 10, the value of objective function remains unchanged. Certainly, under this situation the optimal discharging flow pattern has changed. This means that there are other constraints that control the final objective value.

Surely, the demand plays a crucial role in the model. For example, if we reduce the time-invariant demand rate of ramp 4 from 4 to 3, the objective value will reduce 2531-2449=82. In order to further consider the influence of demand, a different type of demand in a crude form will be considered in subsection 6.3.

*6.2.2. Fixed Total demand of the S-shaped curve*

In this subsection, we will consider the mobilization curve type of demand. Suppose that the total demands are 30, 35, 40, 35, and 30 related to on-ramps 1, 2, 3, 4, and 5, respectively. The evacuation time horizon is still 20. But the effective arrival time limit has changed to 14. Parameters *LR* and *HF* in (17) are set to 0.5 and 5, respectively. Note that choosing proper value for *LR* will be very important. Through calculation, we know that when *LR* is set to 0.6, 0.7, or 0.8, there is no feasible solution to the model. The reason lies in the further concentrated demand in several time intervals. When *LR* is set to 0.3, the clearance time will extended to 15 due to the considerable number of later coming evacuees.

In Table 5, the time dependent demands generated from eq. (17) with *LR*=0.5 and *HF*=5 is given. Note that the truncated data corresponding to time interval 1 are the cumulative demands until the time interval, and the final effective arrival time interval has the residual demand excluding the arrived. Fig. 11 is the corresponding cumulative flow curve. At the beginning of evacuation, the cumulative released flow decreases first, and then increases to the bottleneck capacity. This change corresponds to the specific time dependent demand pattern given in Table 5. In Table 6, the specific released flow pattern w.r.t. equal weights and S-shaped curve demand is also listed.

In figures 12 and 13, we show how the changing weights impact the final released flow pattern. In Fig. 12, we set the weight of on-ramp $r$ to $r^3$ so as to make the upstream ramps with higher priorities of releasing flow. In contrast to Fig. 12, Fig. 13 shows the cumulative flow pattern when we set the weight of on-ramp $r$ to $(9-r)^3$.

It is obvious that the higher weight leads to earlier releasing of more vehicles from the corresponding on-ramp.

**Table 5**

The time dependent demand resulted from S-shaped cure with *LR*=0.5 and *HF*=5.

| Time | 1 | 2 | 3 | 4 | 5 | 6 | 7 | 8 | 9 | 10 | 11 | 12 | 13 | 14 |
|---|---|---|---|---|---|---|---|---|---|---|---|---|---|---|
| Ramp 1 | 2.28 | 1.30 | 1.90 | 2.60 | 3.26 | 3.67 | 3.67 | 3.26 | 2.60 | 1.90 | 1.30 | 0.85 | 0.54 | 0.88 |
| Ramp 2 | 2.66 | 1.52 | 2.21 | 3.03 | 3.80 | 4.29 | 4.29 | 3.80 | 3.03 | 2.21 | 1.52 | 1.00 | 0.63 | 1.03 |
| Ramp 3 | 3.03 | 1.73 | 2.53 | 3.46 | 4.34 | 4.90 | 4.90 | 4.34 | 3.46 | 2.53 | 1.73 | 1.14 | 0.72 | 1.17 |
| Ramp 4 | 2.66 | 1.52 | 2.21 | 3.03 | 3.80 | 4.29 | 4.29 | 3.80 | 3.03 | 2.21 | 1.52 | 1.00 | 0.63 | 1.03 |
| Ramp 5 | 2.28 | 1.30 | 1.90 | 2.60 | 3.26 | 3.67 | 3.67 | 3.26 | 2.60 | 1.90 | 1.30 | 0.85 | 0.54 | 0.88 |

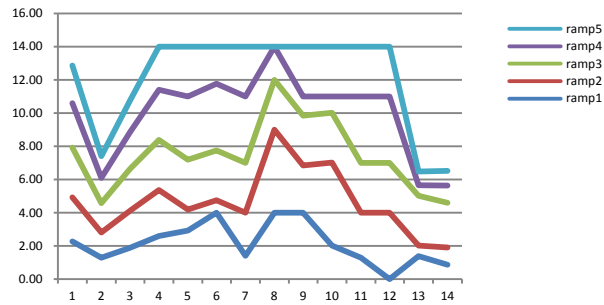

Fig.11. The cumulative flow with equal weights to ramps.

**Table 6**

The released flow pattern w.r.t. equal weights and S-shaped curve demand.

| Time | 1 | 2 | 3 | 4 | 5 | 6 | 7 | 8 | 9 | 10 | 11 | 12 | 13 | 14 |
|---|---|---|---|---|---|---|---|---|---|---|---|---|---|---|
| Ramp 1 | 2.28 | 1.30 | 1.90 | 2.60 | 2.93 | 4.00 | 1.40 | 4.00 | 4.00 | 2.02 | 1.30 | 0.00 | 1.40 | 0.88 |
| Ramp 2 | 2.66 | 1.52 | 2.21 | 2.78 | 1.27 | 0.76 | 2.60 | 5.00 | 2.85 | 5.00 | 2.70 | 4.00 | 0.63 | 1.03 |
| Ramp 3 | 3.00 | 1.77 | 2.53 | 3.00 | 3.00 | 3.00 | 3.00 | 3.00 | 3.00 | 3.00 | 3.00 | 3.00 | 3.00 | 2.70 |
| Ramp 4 | 2.66 | 1.52 | 2.21 | 3.03 | 3.80 | 4.00 | 4.00 | 2.00 | 1.15 | 0.98 | 4.00 | 4.00 | 0.63 | 1.03 |
| Ramp 5 | 2.28 | 1.30 | 1.90 | 2.60 | 3.00 | 2.24 | 3.00 | 0.00 | 3.00 | 3.00 | 3.00 | 3.00 | 0.81 | 0.88 |

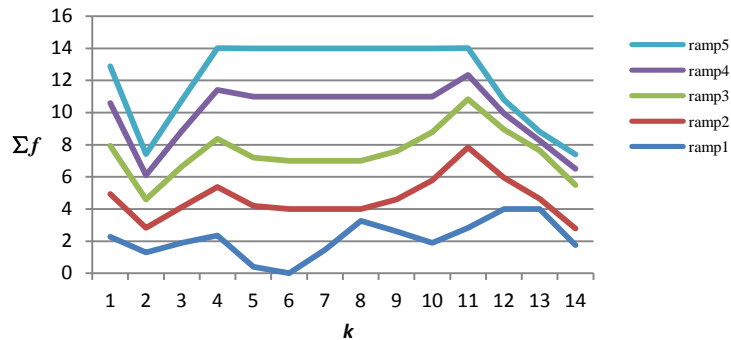

Fig.12. The cumulative flow with weight $r^3$ to ramp $r$.

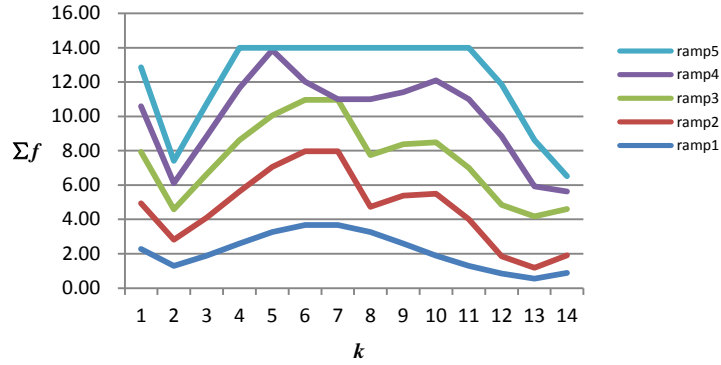

Fig.13. The cumulative flow with equal weight $(9-r)^3$ to ramp $r$.

*6.3. Robust analysis of uncertain demand*

*6.3.1. The uncertainty level and the nominal demand*

In this subsection, we will investigate the impact of uncertain demand in a polyhedral form. The evacuation time horizon is set to 10 and the effective arrival time limit is 5. The total demand upper bounds are 19, 18.5, 18, 21, and 20 corresponding to ramps 1, 2, 3, 4, and 5, respectively. For simplicity of illustration, we assume that all the ramps have the same demand constraints in a polyhedral form. The AARC model (20) has 29,482 variables and 13,057 constraints in this situation.

In Table 7, we give out the objective values with respect to different nominal demands (ND) and uncertainty levels. A figure corresponding to table 7 is drawn below. From Table 7 and Fig. 14, we can see that the higher the uncertainty level, the smaller the objective value becomes if the nominal demand remains the same. With an increasing nominal demand, an improved objective value is obtained. These observed phenomena are reasonable. Since with the increased uncertainty level, it becomes more difficult to supply an optimal solution to model (20). The MAX type of objective function will become smaller under this circumstance. The increased nominal demands means more total demand that will lead to a bigger objective value. From Fig.14, we can also observe the approximate linearity of the change. This will benefit the practical estimate of objective value with respect to the modified nominal demand or uncertainty level. At the end of this subsection, we want to remind readers that the nominal demand has an upper bound w.r.t. a specified uncertainty level due to other constraints. An improper value of nominal demand may lead to non feasible solution to our problem.

**Table 7**

The impacts of uncertainty levels and nominal demands on objective value.

| $\theta$ | Nominal Demand (ND) | | | | | |
|---|---|---|---|---|---|---|
| | 3 | 3.2 | 3.4 | 3.6 | 3.8 | 4 |
| 0.10 | 540 | 570 | 592.5 | 615 | 636 | 654 |
| 0.14 | 516 | 550.4 | 575.5 | 597 | 618.5 | 638 |
| 0.16 | 504 | 537.6 | 567 | 588 | 609 | 630 |
| 0.18 | 492 | 524.8 | 557.6 | 579 | 599.5 | 620 |
| 0.22 | 468 | 499.2 | 530.4 | 561 | 580.5 | 600 |
| 0.26 | 444 | 473.6 | 503.2 | 532.8 | 561.5 | 580 |
| 0.28 | 432 | 460.8 | 489.6 | 518.4 | 547.2 | 570 |
| 0.30 | 420 | 448 | 476 | 504 | 532 | 560 |

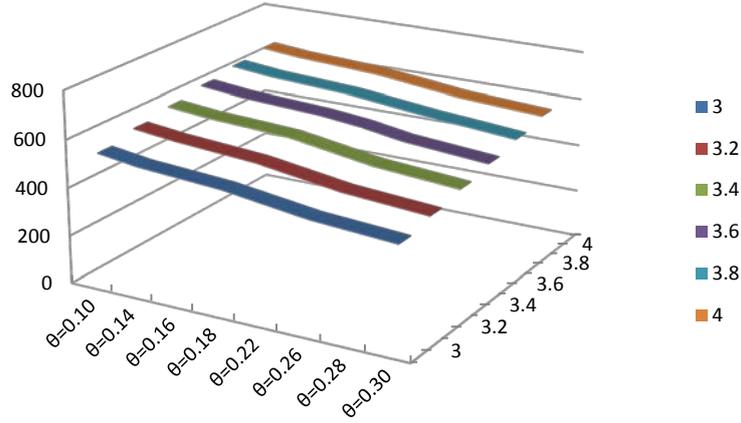

Fig.14. The changing trend of objective value w.r.t. $\theta$ and nominal demand.

*6.3.2. SSP vs. ARRC*

In this subsection, we will compare the sampling based stochastic approach with AARC. The sampling based stochastic programming (SSP) method assumes that a group of demand samples is given. By solving the deterministic model repeatedly w.r.t. all of the samples, the average objective value is obtained as the final estimate of system's performance when using discharging flow pattern w.r.t. the average demand as the evacuation discharging strategy. These average discharging flows will be viewed as the given discharging flows to direct the real evacuation later. This method is simple, but the feasibility of its solution is not guaranteed in many cases. It cannot supply an optimal solution to the worst case.

To compare SSP with AARC, it is necessary to suppose that the demands follow a specified Beta Distribution. When a demand sample is produced using the given Beta Distribution, we first judge if it satisfies the upper bound constraints, and then use them generate a discharging flow pattern according to the following rules. The specific process to carry out a simulation using SSP is described as follows. If the sample produced violates the upper bound constraints, we just drop it and try to generate a new one. With the feasible sample satisfying the polyhedral demand constraints, we can step by step compare the current realized demand to an on-ramp with its given discharging flow pattern and then determine the specified discharging flow for every time interval. If the current demand is greater than the given discharging flow in the time interval in question, the extra part of the demand will be left back to be added into the realized demand of the next time interval. If the current demand is less than the given flow, it will be released totally. If at the end of the evacuation there are some vehicles left back, the corresponding demand will be weighed by a penalty cost (PC), and then be added into the final objective value.

At the beginning of comparison, 50 samples are generated using the specified beta distribution function. They will be used as the sampling demands to produce the estimate of system performance and the given discharging flows. Then 1000 samples are generated to simulate the real evacuation situation. The statistic data are summarized in Tables 8 and 9 corresponding to Beta(1,1) and Beta(2,5), respectively. In the comparison of SSP with AARC, different uncertainty levels are used. In table 8, the common nominal demand is set to 3.6. The penalty cost takes value -5. As is well known, Beta(1,1) is in fact the uniform distribution. In Table 9, the sampling demand is supposed following Beta(2,5). The nominal demand is changed to 3.2. The penalty cost becomes -20.

**Table 8**

AARC vs. SSP when $\theta$ changed. (Beta(1,1), ND=3.6, PC= -5.)

| $\theta$ | Obj | | Avg | | sd | | Worst | |
|---|---|---|---|---|---|---|---|---|
| | SSP | AARC | SSP | AARC | SSP | AARC | SSP | AARC |
| 0.1 | 651.42 | 615 | 641.87 | 646.58 | 3.88 | 3.57 | 627.63 | 635.47 |
| 0.15 | 650.25 | 592.5 | 635.71 | 634.04 | 5.79 | 4.46 | 615.32 | 618.42 |
| 0.2 | 645.90 | 570 | 629.82 | 640.28 | 7.84 | 7.12 | 604.72 | 619.74 |
| 0.25 | 644.12 | 540 | 622.31 | 626.47 | 9.98 | 9.57 | 588.43 | 597.99 |
| 0.3 | 642.20 | 504 | 616.37 | 622.62 | 11.93 | 11.62 | 573.97 | 585.55 |

**Table 9**

AARC vs. SSP when $\theta$ changed.(Beta(2,5) , ND= 3.2, PC= -20.)

| $\theta$ | Obj | | Avg | | sd | | Worst | |
|---|---|---|---|---|---|---|---|---|
| | SSP | AARC | SSP | AARC | SSP | AARC | SSP | AARC |
| 0.1 | 592.71 | 570 | 577.26 | 591.33 | 5.95 | 2.35 | 542.89 | 584.52 |
| 0.15 | 585.26 | 544 | 557.63 | 562.44 | 9.01 | 3.60 | 513.19 | 551.29 |
| 0.2 | 575.60 | 512 | 541.42 | 554.86 | 12.10 | 5.29 | 491.57 | 540.86 |
| 0.25 | 566.16 | 480 | 523.57 | 540.35 | 14.45 | 6.45 | 472.24 | 522.33 |
| 0.3 | 555.61 | 448 | 501.28 | 512.94 | 18.11 | 7.38 | 433.39 | 490.50 |

From the data in Tables 7 and 8, we can see the objective value of SSP is higher than AARC. This is very reasonable because SSP approach only finds the average of maximum objective value with the given samples, but AARC approach generates the best solution to the worst case. When we consider the average and the standard deviation of simulated objective values, we can see AARC approach outperforms SSP approach most of time with a relative small oscillation range. Regarding to the worst case, AARC approach always supplies a better solution than SSP. At last, we want to remind readers that the AARC solution guarantees the feasibility and provides a guaranteed lower bound on the targeted objective value, but the SSP solution guarantees neither of them.

## 7. Conclusions

In this paper, we proposed a linear programming model to deal with a freeway used under an emergency evacuation situation. It is a trade-off between an effective mathematical model and the complex reality. To answer what we mainly concerned in a concise way, we have overlooked many factors, e.g. complex traffic flow states, to formulate a relatively simple LP. This can be viewed as a trial to bridge the gap between the existing complex models for the ordinary situation and the requirement for an effective simple model for an emergency evacuation using a freeway. The constructed models in this paper are made use of to address three issues that are important for an effective emergency evacuation.

The first issue is about realizing the predetermined priorities assigned to some appointed on-ramps. Using sensitivity analysis technique of LP, we can achieve our goal by adjusting the objective function coefficients (or weights) to the predetermined priorities. Four types of predetermined priority patterns have been investigated in Section 3. These analyses will bring flexibility to evacuation managers about setting proper priorities and realizing them in a convenient way.

The second issue we have focused on is about utilization of variable mainline flow capacity. We addressed how to identify the bottleneck of mainline and the possible influence on the system performance due to the changed capacity on the bottleneck. We also suggest a way to roughly estimate the residual capacity on some segments of mainline. These segments may be used to buffer more vehicles and then lighten the pressure of upstream queuing in on-ramps.

Last, but not the least, we address the issue about demand uncertainty. In this part, we assume that a polyhedral form of demand constraints is adopted. By utilizing the technique of AARC, we have transformed such a semi-infinite intractable model into a tractable LP. Through comparison with sampling based stochastic methods, the advantage of this formulation is clarified. The application of the model in practice is promising.

To make the methods proposed in this paper more practical, further study in several directions is required. The field test is welcome to further convince related models' effectiveness and accuracy. Some simulation experiments are desired to adjust these simple models to complex nonlinear situations. How to make use of the residual capacity in mainline is an interesting issue desiring further efforts.